\documentclass[11pt,leqno]{amsart}
\usepackage{amsfonts,amssymb}
\usepackage{graphicx,tikz}
\usepackage{tikz-cd} 
\usepackage{dsfont}
\usepackage{float}
\usepackage{amsmath, amsthm, amsfonts}
\usepackage{hyperref}
\usepackage{enumitem}  
\usepackage[all,cmtip]{xy}
\usepackage[capitalise]{cleveref}
\usepackage{comment}
\usepackage{enumitem}
\usepackage{amsrefs}
\usepackage{nicefrac}
\usepackage[official]{eurosym}
\usepackage{graphicx}
\usepackage{nomencl}
\usepackage{amsfonts}
\usepackage{hyperref}
\usepackage{tabto}
\usepackage{amssymb}
\usepackage{fancyhdr}
\usepackage{amscd}
\usepackage[english]{babel}

\usepackage[utf8]{inputenc}

\usepackage{hyperref}
\usepackage[normalem]{ulem}
\usepackage{todonotes, comment,graphicx}

\usepackage{amssymb,amsmath,amsthm, textcomp}
\theoremstyle{plain}
\newtheorem{theorem}{Theorem}

\newcommand{\problem}[2]{\theoremstyle{plain}
   \newtheorem*{thm#1}{#1}\begin{thm#1}#2\end{thm#1}}

\theoremstyle{remark}

\newtheorem*{rem*}{Remark}

\theoremstyle{definition}
\newtheorem{definition}{Definition}[section]

\def\a{\alpha}

\def\ba{\breve{\alpha}}
\def\be{\breve{\eta}}

\def\ha{\hat{\alpha}}
\def\hb{\hat{\beta}}

\title[Simple Groups and Cryptography]{Applications of Finite non-Abelian Simple Groups to Cryptography in the Quantum Era}
\author[M. González Vasco]{María Isabel González Vasco}
\author[D. Kahrobaei] {Delaram Kahrobaei} 
\author[E. McKemmie]{Eilidh McKemmie}

\address{Catedratica de Universidad 
Departamento de Matematicas, Universidad Carlos III de Madrid 
Campus de Legan\'es, Madrid, Spain}

 \address{Departments of Computer Science and Mathematics, Queens College, City University of New York, USA, 
     Department of Computer Science, University of York, UK, 
     Initiative for the Theoretical Sciences, Graduate Center,  City University of New York, USA, Department of Computer Science and Engineering, Tandon School of Engineering, New York University, USA}
\address{Mathematics Department, Rutgers University, \\New Brunswick, New Jersey, USA}

\begin{document}

\maketitle
\begin{abstract}


  The theory of finite simple groups is a (rather unexplored) area likely to provide interesting computational problems and modelling tools useful in a cryptographic context. In this note, we review some applications of finite non-abelian simple groups to cryptography  and discuss different scenarios in which this theory is clearly central, providing the relevant definitions to make the material accessible to both cryptographers and group theorists, in the hope of stimulating further interaction between these two (non-disjoint) communities. In particular, we look at constructions based on various group-theoretic factorization problems, review group theoretical hash functions, and discuss fully homomorphic encryption using simple groups. The Hidden Subgroup Problem is also briefly discussed in this context.
\end{abstract}
\section{Introduction}

Cryptography is built upon the computational hardness of certain mathematical problems. One of the main tools within this area are \emph{one-way} functions (informally, functions that can be efficiently evaluated while there are no efficient methods to compute preimages, possibly unless there is a secret key giving additional information). Computational tasks like factoring large integers or decoding with respect to  random codes are flagship examples of mathematical problems naturally defining one-way functions. Of course, considering different computational models has a large impact in how such cryptographic-amenable problems can be selected; in particular, since the 1980s the appearance of quantum computing has necessitated the search for problems that will remain hard even if a quantum computer is available. The field of \emph{post-quantum} cryptography revolves around cryptographic designs whose security relies on these kind of problems.

There have been many cryptographic proposals based on problems in group theory, see the recent book and survey by Kahrobaei et al \cite{KFNHB, KFN}. While it is not easy to classify problems as quantum resistant in a reasonable way, we do know of some problems that quantum computers can tackle with a significant advantage. The main menance is Shor's~\cite{Shor} quantum algorithm, which gives an exponential gain for solving problems that fit a certain ``period-finding'' description. Factoring large integers or solving discrete logarithms in finite cyclic groups fall into this category. Remarkably, it seems that the ideas behind Shor's algorithm can be extended to exploit normal subgroup structure in other groups. Simple groups are those with no non-trivial normal subgroups, so it is natural to ask whether finite simple groups may be harder than other groups for quantum computers to deal with. This leads us to suggest that the finite simple groups may be a good setting for post-quantum cryptographic schemes.

In the literature there are proposals using finite non-abelian simple groups for constructing many different tools: encryption and digital signature schemes, fully homomorphic encryption designs and hash functions. In this survey, we will take a closer look at the status of some proposed applications of the theory of finite simple groups to the design of hash functions, public-key encryption and fully homomorphic encryption. Our aim is not to be exhaustive but simply to give the reader a glimpse of the vast amount of unexplored avenues within this area, with a focus on some challenging group-theoretic and computational problems relevant to building sound cryptographic constructions.



\medskip

\noindent{\emph{Paper Roadmap.}} We start with a brief introduction to the finite simple groups and their classification in Section~\ref{sect: CFSG}. In Section~\ref{sect: Cayley}, we introduce Cayley hash functions and give an example cryptographic construction. We then discuss the difficulty of a certain factorization problem in groups that is linked to their security, and a related group theoretic conjecture. In Section~\ref{sect: LS}, we define logarithmic signatures and another factorization problem in groups which has been used as justification for several public key cryptosystems. We give an example of a cryptographic construction and discuss a related group-theoretic conjecture. Section~\ref{sect: FHE} discusses fully homomorphic encryption schemes and a method of building them from homomorphic encryption on groups, while in Section~\ref{sect: HSP} we discuss the Hidden Subgroup Problem for cryptanalysis of proposed schemes using finite non-abelian simple groups against possible quantum attacks. Section~\ref{sect: open problems} concludes the paper with a summary of the exciting open problems we discussed.


\section{Preliminaries: Finite Simple Groups}\label{sect: CFSG}

A \emph{simple group} is a non-trivial group whose only normal subgroups are itself and the trivial group. We are also interested in some \emph{quasisimple} groups: $G$ is quasisimple if it is perfect (i.e. equal to its own commutator subgroup $G=[G,G]$) and its group of inner automorphisms ${\rm Inn}(G)$ is simple. We focus here on finite groups since our cryptography applications require finite data structures.

There is a classification of all finite simple groups whose proof was completed in the 2000s after many years of work by a large number of mathematicians. For a brief historical overview, see \cite{Aschbacher2004}. The list of finite simple groups is as follows:

\begin{theorem}
If $G$ is a finite simple group then either $G$ is abelian, in which case it is a cyclic group of prime order, or $G$ is non-abelian, in which case one of the following holds:
\begin{itemize}
    \item $G\cong A_n$ is an alternating group on $n\ge 5$ letters
    \item $G$ is a group of Lie type
    \item $G$ is one of $26$ sporadic groups.
\end{itemize}
\end{theorem}
The proof takes up many books, see for example the series \cite{GLS}. For a more introductory textbook describing all the groups in detail, see \cite{wilson2009finite}. 



The groups of Lie type are the classical groups and the exceptional groups over finite fields. We describe these groups briefly here, and refer the reader to a standard textbook by Carter \cite{C1989} for more details. These groups are defined over finite fields. We use $p$ to denote the characteristic of the field, which is a prime, and $q$ to denote the order of the field, which is a power of $p$. Each finite group of Lie type has an underlying root system which determines an integer known as the \emph{rank} of the group.

The classical groups are those which are natural matrix groups, and there are four types for every integer $n\ge 2$ and prime power $q$. For example, the projective special linear group of $n\times n$ matrices over a field of order $q$, denoted $PSL_n(q)$, has rank $n-1$ and is simple except when $n=2$ and $q=2, 3$. The other classical groups are the groups of unitary, orthogonal and symplectic matrices over finite fields. 
We are also interested in finite quasisimple classical groups, for example the special linear group $SL_n(q)$. In characteristic $2$ we have that $SL_n(2^k)=PSL_n(2^k)$ which is simple for $k>1$.

The exceptional groups do not have such natural representations as groups of matrices, and all have rank at most $8$. There are 10 infinite families indexed by prime powers $q$. One such family is the Suzuki groups which are defined over fields of order $2^{2n+1}$ which we denote by $Sz(2^{2n+1})$.


\section{Factorization Problem and Cayley Hash Functions}\label{sect: Cayley}
A hash function is a function whose input is an arbitrarily large message and whose output is a fixed-length \emph{hash}. Hash functions are a cryptographic primitive with a variety of cryptographic applications, each requiring different security properties (see any cryptography textbook, for example \cite[Chapter~6]{Katz2021}). Desirable properties of a hash function $h:M \rightarrow N$ include \emph{preimage resistance} -- given $n \in N$ it should be computationally infeasible to find $m \in M$ such that $h(m)=n$ -- and \emph{collision resistance} -- it should be computationally infeasible to find $m \ne m' \in M$ such that $h(m)=h(m')$.

Z{\'{e}}mor \cite{Zemor} defined group theoretic hash functions based on Cayley graphs of finitely-generated groups, following work of Bosset and Camion \cite{Camion_1986}.

\begin{definition}
Let $G$ be a finitely generated group with a generating set $S=\{g_1, ..., g_k\}$ which is closed under taking inverses.

\begin{itemize}
    \item The Cayley graph $\Gamma(G,S)$ is a graph with vertex set $G$ and an edge from $g$ to $h$ if and only if $g=g_i h$ for some $i$.
    \item The Cayley hash function $h_{G,S}:\{1, ..., k\}^* \rightarrow G$ is defined by \\ $h_{G,S}(m_1, m_2, ..., m_r)=g_{m_1}g_{m_2}\cdots g_{m_r}$. We refer to $(m_1, m_2, ..., m_r)\in \{1, ..., k\}^*$ as a \emph{word of length} $r$.
\end{itemize}
\end{definition}

Note that evaluation of $h_{G,S}$ at $(m_1, m_2, ..., m_r)$ corresponds to traversing the path $(1, g_{m_1}, g_{m_1}g_{m_2}, \cdots, g_{m_1}g_{m_2}\cdots g_{m_r})$ in the Cayley graph $\Gamma(G,S)$.

Preimage resistance for Cayley hash functions is equivalent to the difficulty of writing a given element of $G$ as a product of elements of $S$, or finding a path from $1$ to the given element in the Cayley graph. This is called the Factorization Problem.

\problem{Factorization Problem}{Given $h \in G$ find a ``short'' word $(m_i)_i$ such that $\prod_i g_{m_i}=h$. Equivalently, given $h \in G$ find a ``short'' path from $1$ to $h$ in $\Gamma(G,S)$.}

It should be noted that finding minimal such words or paths is the NP-hard Minimum Generator Sequence Problem \cite{EVEN1981311}.

\subsection{Cryptographic constructions}

There have been several choices of generating sets proposed for Cayley hash functions over $SL_2(q)$ \cite{Zemor, Tillich, Zemor_1994,Petit_2013, Petit2016}, but there are known attacks in each case \cite{Petit_2008, Grassl_2010, Petit_2011, tinani2023methods}. Recently, Le Coz, Battarbee, Flores, Koberda and Kahrobaei \cite{LBFKK} proposed a generating set for the quasisimple group $SL_n(p)$ for prime $p$. The Factorization Problem in this case can be reduced to solving a system of $n^2$ multivariate polynomial equations in $O(\log p)$ unknowns over $\mathbb{F}_p$ \cite[Section~3.2]{LBFKK} which is known to be NP-hard in the worst case.

As an example, we describe a particularly simple scheme proposed by Z{\'{e}}mor \cite{Zemor}.

Let $p$ be a prime and associate to the bit $0$ the matrix $A=\begin{pmatrix}1&1\\0&1\end{pmatrix} \in SL_2(p)$ and to the bit $1$ the matrix $B=\begin{pmatrix}1&0\\1&1\end{pmatrix}\in SL_2(p)$. Then the hash function $h_{SL_2(p), \{A, B\}}$ sends a binary number of arbitrary length to the appropriate product of $A$s and $B$s.

These parameters were chosen to allow efficient evaluation of the hash function, but the resulting hash function is not collision resistant:  Tillich and Z{\'{e}}mor \cite{Tillich1994,Tillich} show it is possible to find many factorizations of the group identity. Inserting any such factorization into any word gives a collision.

\subsection{Progress towards solving the Factorization Problem}\label{sect: Cayley proofs}

Babai and Seress conjectured \cite{Babai1992} that short paths exist in the Cayley graphs of finite simple groups:

\problem{Babai's conjecture}
{There exists a constant $c > 0$ such that, for any $h$ in a finite simple non-abelian group $G$, and any generating set $S$, there is a path from $1$ to $h$ in $\Gamma(G,S)$ of length at most $(\log |G|)^c$. That is, every element of $G$ may be written as a word of length at most $(\log |G|)^c$ in the elements of $S$.}

For groups of Lie type of bounded rank, Babai's conjecture has been proved by Helfgott, Pyber, Szabo, Breuillard, Green and Tao \cite{Helfgott2008, Pyber_2014, Breuillard_2011}. The remaining cases are the alternating groups (for which Helfgott and Seress \cite{Helfgott_2014} have the best bound) and groups of Lie type of unbounded rank. In many cases there are partial results proving Babai's conjecture for certain generating sets. For example Babai and Hayes \cite{Babai2005} prove Babai's conjecture for almost all generating sets of alternating groups, and Eberhard and Jezernik recently showed \cite{Eberhard_2021} that Babai's conjecture holds for large rank groups of Lie type for almost all large enough sets $S$. See \cite[Section~1]{Eberhard_2021} for more details on the current status of Babai's conjecture.

Babai's conjecture would imply that for every $h \in G$ there is a path of length $(\log |G|)^{O(1)}$ from $1$ to $h$ in the Cayley graph, and the goal of cryptanalysts is to explicitly construct such short paths, while the goal of cryptographers is to find generating sets that make this as difficult as possible. There has been much activity in this area: Minkwitz \cite{Minkwitz_1998} provided an optimization for the Schreier-Sims algorithm \cite{SIMS1970169, Knuth_1991} for solving the Factorization Problem in permutation groups. Babai and Hayes \cite{Babai2005} (see also \cite{Babai2004}) give a Las Vegas algorithm based on a random walk which is able to factorize elements of $A_n$ for almost all generating sets, and Kalka, Teicher and Tsaban \cite[Section~5]{Kalka_2012} provide an algorithm which conjecturally and experimentally gives even shorter words. Babai, Kantor and Lubotzky \cite{Babai_1989} showed that every finite simple non-abelian group $G$ has a set of generators $S$ of size at most $7$ for which there is an algorithm that finds words of length $O(\log |G|)$ in $O(\log |G|)$ time. Of the groups of Lie type, $PSL_n(q)$ and $SL_n(q)$ have been most closely studied and there are a handful of specially chosen generating sets for which there are efficient algorithms \cite{Larsen_2003, Riley2005, Kassabov_2007, Petit_2008}. Another approach of Kantor and Seress and Dietrich, Leedham-Green and O'Brien is to represent classical groups as so-called black-box groups 
and use a Las Vegas algorithm to attempt to construct standard generating sets in which to solve the Factorization Problem \cite{kantor2001black, Dietrich_2015}. For all generating sets of $SL_2(2^k)$ there is a subexponential-time algorithm giving subexponential-length words \cite{Petit_2012}. However, there is no efficient algorithm which works for all groups and generating sets.

\section{Public Key Constructions from Logarithmic Signatures}\label{sect: LS}


Since the 1980s, there have been several attempts to exploit the computational properties of so-called factorization 
sequences of finite groups to derive one-way functions, including trapdoor functions -- one-way functions for which it becomes easy to compute preimages given some extra information (see for instance~\cite{DBLP:phd/dnb/Reichl15a}).

\begin{definition}
Let $G$ be a finite group. We may identify $G$ with a permutation group acting on $n$ points where $n\le |G|$. Call this $n$ the \emph{degree} of $G$. Fix $s \in \mathbb{N}$ and for each $i=1, ..., s$ let $\alpha_{ij}\in G$ and consider $\alpha = (\alpha_1,\dots, \alpha_s)$ where $\alpha_i=(\alpha_{i1},\dots, \alpha_{in_i})$. We denote by $\ell(\alpha)=\sum_{i=1}^sn_i$ the \emph{length} of $\alpha$.

We say that $(i_1,\dots, i_s)\in \mathbb{N}^s$ is a \emph{factorization sequence for $g\in G$ w.r.t. $\alpha$} if $g=\alpha_{1i_1}\cdots \alpha_{si_s}.$ Denote by $n[\alpha,g]$ the number of different factorization sequences for $g$ induced by $\alpha.$ We say that $\alpha$ is a
\begin{itemize}
\item \emph{cover} if $n[\alpha, g] >0$ for any $g\in G$.
\item \emph{logarithmic signature} if $n[\alpha, g] =1$ for any $g\in G$. A logarithmic signature $\alpha$ is called \emph{tame} if factorization sequences may be computed in polynomial time in the degree of $G$ for every $g$ w.r.t. $\alpha$, and \emph{wild} otherwise.
\end{itemize}
Note that by definition $\alpha$ is a logarithmic signature if and only if $\alpha$ is a cover and $\prod_{i=1}^s n_i = |G|$.
\end{definition}

If the group law can be computed efficiently, it is ``easy" to construct group elements by simply selecting one element from each $\alpha_i$; the reverse process may, however, be rather involved computationally. The next section reviews several proposals exploiting this dichotomy to define useful one-way functions.

\subsection{Cryptographic Constructions}

The first private-key cryptographic construction using factorization sequences was PGM (Permutation Group Mappings) which was proposed by Magliveras~\cite{DBLP:journals/joc/MagliverasM92} and uses logarithmic signatures for permutation groups to create one-way functions. Later, Magliveras et al.~\cite{MaStTr00} proposed ${\rm MST_1}$, a public-key cryptosystem built upon the same idea with an additional trapdoor for the one-way function of PGM. They also proposed a variant called ${\rm MST_2}$ based on a special kind of cover called a \emph{mesh}. 
Later, Lempken et al.~\cite{DBLP:journals/joc/LempkenTMW09} proposed
${\rm MST_3}$ based on the difficulty of factoring group elements with respect to random covers for large subsets of finite non-abelian groups with large center.

We now give a description of $MST_1$, the simplest of these constructions. For a natural number $m$ we denote by ${\mathbb Z}_m=\{0,1,\dots,m-1\}$ the ring of integers modulo $m$. Fix a finite permutation group $G$ and a tame logarithmic signature $\eta$ for $G$, both publicly known. For any logarithmic signature
$\alpha=(\alpha_1,\dots,\alpha_s)$ we construct the mappings
$$\begin{array}{rccc}
\lambda: &{\mathbb Z}_{n_1}\times\dots\times {\mathbb Z}_{n_s}&\longrightarrow& {\mathbb Z}_{|G|}\\
         &(r_1,\dots,r_s)&\longmapsto&\sum_{i=1}^s\left(r_i\cdot\prod_{j=1}^{i-1
}n_j\right)
\end{array}
$$

and
$$\begin{array}{rccc}
\Theta_{\a}: &{\mathbb Z}_{n_1}\times\dots\times {\mathbb Z}_{n_s}&\longrightarrow&G\\
             &(r_1,\dots,r_s)&\longmapsto&\alpha_{1 r_1}\cdots\alpha_{s r_s}
\end{array},$$
which one may check are bijective. Thus, the functional
composition of $\Theta_{\a}$ and $\lambda^{-1}$ yields a bijection
\[
\begin{array}{rccc}
\ba:&
\mathbb{Z}_{|G|}&\longrightarrow &G\\
&n&\longmapsto&(\Theta_{\a}\lambda^{-1})(n)=\Theta_{\a}\left(\lambda^{-1}(n)
\right).
\end{array}
\]We will use $\be^{-1}$ to identify $G$ with $\mathbb{Z}_{|G|}$, allowing us to associate to each 
logarithmic signature $\alpha$ a permutation $\ha:=\be^{-1}\ba\in S_{|G|}$.

\begin{figure*}[h]
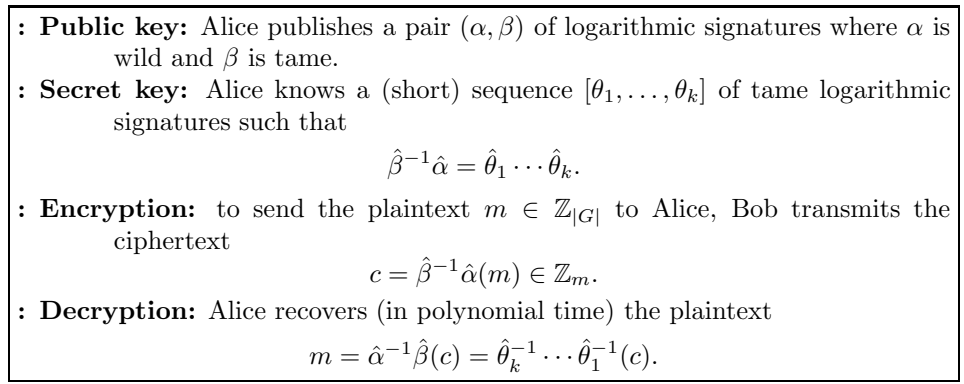

\fbox{\begin{minipage}{0.98\textwidth}
\small
\begin{description}
\item{\bf Public key:} Alice publishes a pair
$(\alpha,\beta)$ of logarithmic signatures where $\alpha$ is wild and $\beta$ is tame.

\item{\bf Secret key:} Alice knows a (short) sequence
$[\theta_1,\dots,\theta_k]$ of tame logarithmic signatures such that
$$\hb^{-1}\ha=\hat\theta_1\cdots\hat\theta_k.$$

\item{\bf Encryption:} to send the plaintext $m\in{\mathbb
Z}_{|G|}$ to Alice, Bob transmits the ciphertext
$$c=\hb^{-1}\ha(m)\in{\mathbb Z}_m.$$

\item{\bf Decryption:} Alice recovers (in polynomial time) the
plaintext
$$m=\ha^{-1}\hb(c)=\hat\theta_k^{-1}\cdots\hat\theta_1^{-1}(c).$$
\end{description}
\end{minipage}}
\caption{$MST_1$ encryption scheme}
\label{Fig:MST1}
\end{figure*}

For ${\rm MST_1}$, the public key is a wild logarithmic signature $\alpha = (\alpha_1,\dots, \alpha_s)$ and a tame logarithmic signature $\beta = (\beta_1,\dots, \beta_s)$ for the same group $G$. The private key consists of a sequence $[\theta_1,\dots,\theta_k]$ of tame logarithmic signatures such that $\hb^{-1}\ha=\hat\theta_1\cdots\hat\theta_k$, which opens the trapdoor to efficient computation of factorization sequences w.r.t. $\alpha$. As discussed in~\cite{DBLP:journals/joc/MagliverasM92}, it is not known how to efficiently compute an appropriate sequence $[\theta_1,\dots,\theta_k]$. 
The encryption scheme is depicted in Figure~\ref{Fig:MST1}.

\subsection{Producing hard factorizations}\label{sect: LS constructions}
All the above constructions base their security on the claimed hardness of computing factorizations of group elements with respect to some public cover. To support such a claim,  the problem of factoring w.r.t a cover should be reduced as closely as possible to another computational problem that we can ``safely" assume to be hard enough.

In the construction of ${\rm MST_1},$ a critical point is the choice of the public wild logarithmic
signature $\alpha$ along with a trapdoor (the factorization into the tame logarithmic signatures $\theta_i$ ($1\le i\le k$)). Magliveras et al.~\cite{MaStTr00} suggested picking $\alpha$ to be a \emph{totally-non-transversal} logarithmic signature, meaning that none of the $\alpha_i$ is a coset of a non-trivial subgroup of $G$. This was later proven in~\cite{BGMVS} to be insufficient since for $n\ge 5$ there are tame totally-non-transversal logarithmic signatures for all alternating groups $A_n$ and symmetric groups $S_n$.

Similarly, the security of $MST_3$ was questioned in~\cite{DBLP:journals/dcc/VascoPD10} and further cryptanalyzed in~\cite{DBLP:journals/iacr/BlackburnCM09}, where it was proven that factoring with respect to the random covers used is not always a hard problem. While further schemes have been proposed in recent years (see, for instance, \cite{8907845,DBLP:journals/jmc/SvabaT10}) at the writing of this survey, we are unfortunately not aware of a secure method for inducing hard group factorizations suited for cryptographic purposes.


\subsection{In search of Minimal Length Logarithmic Signatures}\label{sect: MLS proofs}

Cryptographic applications motivate nice group-theoretic questions. For example, since the length of covers is a relevant parameter in real-life implementations, one may ask what the minimal length of a logarithmic signature can be, and try to construct logarithmic signatures of this length.

Let $G$ be a finite group of order $|G| = \prod_{j=1}^kp_j^{a_j}$ with $p_1,\dots,p_k$ distinct primes. Gonz{\'a}lez Vasco and Steinwandt~\cite{GVSt01} showed that for each logarithmic signature $\alpha$ for $G$ we have
 \begin{equation}\label{equ:thebound}
    \ell(\alpha) \ge \sum_{j=1}^ka_jp_j,
 \end{equation} and defined a \emph{minimal length logarithmic signature} $\alpha$ to be a logarithmic signature for which equality in (\ref{equ:thebound}) holds. Then they constructed minimal length logarithmic sequences for symmetric and solvable groups. It is not yet known if minimal length logarithmic signatures exist for each finite group, although Magliveras~\cite{Mag01} reduced the problem to simple groups, showing that a minimal counterexample of a group without a minimal length logarithmic signature  must be simple. He also constructed minimal length logarithmic signatures for the alternating groups. The work in \cite{GVSt01, Mag01} leads to the following conjecture for which a constructive proof is desired.

\problem{MLS Conjecture}{
Every finite simple group has a minimal length logarithmic signature.
}


This conjecture remains open in general, but has been proved in several cases. The constructive proofs for symmetric and alternating groups are in essence obtained by the same technique: given a permutation representation of a group $G$, identify a point $P$ so that its stabilizer $G_{P}$ can be factored through a minimal length logarithmic signature and such that there exists a complete set of representatives of $G$ modulo $G_{P}$ which moves $P$ cyclically. The underlying idea is to factor the group into a `product of disjoint pieces' for which a minimal length logarithmic signature exists. In the case that these `disjoint pieces' are two subgroups, this is a rewriting of the group as a {\em knit (or Zappa-Sz\'ep) product}~\cite{Mi,DBLP:journals/em/VascoRS03}.

Lempken and van Trung~\cite{DBLP:journals/em/LempkenT05} use \emph{double coset decomposition} to find minimal length logarithmic signatures for a number of special linear groups and projective special linear groups. Constructions of minimal length logarithmic sequences for all of the simple linear and symplectic groups, as well as some orthogonal groups, are found in~\cite{DBLP:journals/dcc/SinghiS11,DBLP:journals/dcc/SinghiSM10}. These papers consider the action of the group on the natural module, looking at point stabilizers and geometric objects called \emph{spreads}.
Furthemore, Holmes~\cite{Holmes} produced minimal logarithmic signatures for the sporadic groups $J_1$, $J_2$, $HS$, $McL$, $He$ and $Co_3$. Rahimipour, Ashrafi and Gholami \cite{Rahimipour_2015, Rahimipour_2016, Rahimipour_2021} treat the cases of the sporadic groups $J_3$, $Fi_{22}$, $Ru$ and $Suz$, as well as the Tits group $^2F_4(2)'$, the Ree groups $^2G_2(3^{2n+1})$, and some unitary and exceptional groups.

\section{Fully Homomorphic Encryption Schemes}\label{sect: FHE}
Broadly, homomorphic encryption enables computation over encrypted data. A \emph{fully homomorphic encryption (FHE)} procedure is an encryption algorithm $E$ taking as input an element from a ring  $(R,+,\cdot)$ and producing an output in another ring  $(S,+,\cdot)$ such that $E(r+s)=E(r)+E(s)$ and $E(r\cdot s)=E(r)\cdot E(s)$. Such an encryption mechanism allows a third party to do any computations involving $+$ and $\cdot$ without ever decrypting the data. For example, one can take the boolean circuit $(\{0,1\},AND, XOR)$ as the ring, so that a fully homomorphic encryption function respects both $AND$ and $OR$.

There are several known encryption schemes  on rings $(\mathbb{Z}_n,+,\cdot)$  which allow homomorphic computation of only one of the two operations, for example textbook RSA, ElGamal and Goldwasser-Micali, but it appears far more difficult to construct a fully homomorphic scheme. For a detailed survey see \cite{WNK}.\

The most widely known existing fully homomorphic
encryption scheme appeared originally in the thesis of Craig Gentry \cite{10.1145/1536414.1536440}. The security of this solution relies on variants of the so called \emph{bounded-distance decoding}
problem. This problem enjoys a very relevant property for cryptographic purpose, namely, it is \emph{random self reducible},   which
basically means that it is about as hard on average as it is in the
worst case. While this property  allows for (practically meaningful) security proofs, it is unfortunaly the case that the resulting homomorphic encryption algorithm is too inefficient to be practical. Very informally, the reason is that, to
provide semantic security, encryption has to be randomized, but on
the other hand, a homomorphism  should   map zero to zero. To
resolve this conflict, the ciphertext zero is ``masked by noise".
The problem now is that during any computation on encrypted data,
this ``noise" tends to accumulate and has to be occasionally reduced
by re-encryption (also known as {\it bootstrapping}), a process that
produces the equivalent ciphertext but with less noise. This 
is an expensive procedure, and its results in real-life computation being prohibitively slow.

The  quest for more efficient techniques to overcome this issue has resulted in a number of rather efficient schemes. For instance,  in \cite{DBLP:journals/siamcomp/BrakerskiV14,DBLP:conf/eurocrypt/GentryHS12}
 a much slower growth of the noise during homomorphic computations was achieved, providing enough efficiency for practical applications.   Later, in 2013, Gentry, Sahai and Waters~\cite{DBLP:conf/crypto/GentrySW13}  put forward the GSW scheme, a new method to derive more efficient FHE schemes. These techniques were further improved to develop efficient ring variants of the GSW scheme \cite{DBLP:conf/asiacrypt/ChillottiGGI16}. New efficient constructions are constantly being proposed (see \cite{DBLP:conf/eurocrypt/LeeMKCDEY23}), and fully homomorphic encryption is indeed a reality in many practical applications.

\subsection{Simple groups and Fully Homomorphic Encryption}\label{sect: FHE simple groups}

The relevance of finite non-abelian simple groups to fully homomorphic encryption is that they open a door to designing new noise-free fully homomorphic encryption schemes, thus with the potential of being much more efficient than those needing some sort of bootstrapping. 

This idea is quantified by the following theorem of Werner \cite{werner1973finite}. 

\begin{theorem}[\cite{werner1973finite, Ostrovsky2007}]
    There is a fully homomorphic encryption scheme (over a non-zero ring) if and only if there is a finite non-abelian simple group over which there is a homomorphic encryption scheme.
\end{theorem}

Ostrovsky and Skeith gave a constructive proof of this theorem \cite[Corollary~4.26]{Ostrovsky2007}, see \cite[Section~6]{Khamsemanan2016} for more discussion. To construct a noise-free fully homomorphic encryption scheme from a group homomorphism $\phi:G \rightarrow H$, Ostrovsky and Skeith pick an element $g \in G$ of order $2$ and identify the bit $0$ with the identity of $G$, and the bit $1$ with the element $g$. Since any binary function can be written as compositions of the $NAND$ function, it is enough to construct $NAND$ in the group. Ostrovsky and Skeith's proof gives a general formula, and they display an example for the group $A_5$. The details for $A_n$ for $n\ge 6$ are especially short, so we describe them here.

Let $g=(1\,2)(3\,4)$ and $e$ be the identity permutation. For $a, b\in \{e, g\}$. We will give a formula for $NAND(a,b)$. We follow Ostrovsky and Skeith's proof, noting that \[g=[(1\,2)(5\,6),(1\,4)(2\,3)]=[g^{(3\,5)(4\,6)},g^{(2\,4)(5\,6)}].\] Therefore \begin{align*}NAND(a,b)&=g [a^{(3\,5)(4\,6)},b^{(2\,4)(5\,6)}]\\&=(1\,2)(3\,6\,4\,5)a(3\,6\,2\,4\,5)b(2\,6\,3\,5\,4)a(3\,6\,2\,4\,5)b(2\,4)(5\,6). \end{align*}

Armknecht, Gagliardoni, Katzenbeisser and Peter \cite{armknecht2014general} give an attack using quantum computers that undermines the security of any homomorphic encryption scheme whose plaintext and ciphertext spaces are abelian groups, thereby showing that it is impossible to have a quantum secure group homomorphic encryption scheme in this scenario. We are not aware of any literature proposing homomorphic encryption over non-abelian groups, but this is a research avenue worth exploring (see \cite{Nuida_2020} for more discussion).


\section{Hidden subgroup problem: post-quantum analysis}

\label{sect: HSP}

The search for quantum-resistant alternatives to today's common public-key constructions is extremely active. As we mentioned in the introduction, it is of paramount importance to identify and understand which mathematical problems are hard enough in a ``post-quantum'' sense. The \emph{Hidden Subgroup Problem} (HSP) is a generic formulation englobing many such potentially hard problems. HSP can be seen as a way to understand the power of quantum algorithms and the limits of Shor's algorithm in group theoretical language.

\problem{Hidden Subgroup Problem (HSP)}{Given a finitely generated group $G$, a finite set $S$ and an efficiently computable function $f:G \rightarrow S$ such that $f$ is constant and distinct on left cosets of a subgroup $H\le G$ of finite index, find a generating set for $H$.}


Famously, Shor's \cite{Shor} polynomial-time quantum algorithms for the Integer Factorization Problem and Discrete Logarithm Problem rely on a polynomial-time quantum algorithm for HSP in finite cyclic groups and groups of the form $\mathbb{Z}_p\times \mathbb{Z}_p$ for prime $p$. There are efficient quantum algorithms for HSP for all finite abelian groups and for a few classes of finite non-abelian groups. We describe some relevant cases here. See \cite{HK2018} for a full survey.

Hallgren, Russell and Ta-Shma \cite[Theorem~2]{Hallgren_2000} gave a quantum algorithm for finding hidden normal subgroups. This result says nothing about finite simple groups since they have no non-trivial normal subgroups. Kuperberg in \cite{kuperberg2005subexponential}, and Regev in \cite{regev2004subexponential} give subexponential-time quantum algorithms for HSP in dihedral groups. Kuperberg's algorithm requires quantum space $2^{O({\log r})}$, while a generalized version of Regev's in \cite[Theorem~5.2]{childs2014constructing} is slower but less space-expensive. In \cite{battarbee2022subexponential}, the authors extend these algorithms to construct a subexponential quantum algorithm for solving the Discrete Logarithm Problem in semi-direct products.

While we have efficient algorithms in some cases, providing solutions for HSP for all finite groups is considered one of the most important challenges in post-quantum cryptography. A solution to HSP in a finite group implies a solution in all subgroups. Since every finite group is a subgroup of a symmetric group, a solution to HSP for all finite groups is equivalent to a solution to HSP for symmetric groups. Note, however, that the representation of our group $G$ as a subgroup of a symmetric group is relevant here, since if the dimension is large (for example if we consider the group $G$ to be in $S_{|G|}$) we will see exponential blow-up in size and parameters. 

     


{Many of the techniques that have been successfully employed in the above-mentioned cases have been shown to fail for symmetric groups \cite{KempeShalev,kuperberg2005subexponential,Moore2007,Moore2008}. See \cite[Section~3.2]{Alagic2017} for more discussion. Often the obstructions are large subgroups and high-dimensional irreducible representations. Therefore, many of the difficulties in the symmetric case also affect the classical group case \cite{Hallgren2010, Moore2007}.


Understanding the complexity of HSP in finite non-abelian groups is a significant open question with strong connections to many well-known hard problems. This suggests study in this area could unearth one-way functions for the design of post-quantum cryptosystems.





\section{The road ahead: some open problems}\label{sect: open problems}

We have presented different problems related to non-abelian finite simple groups. We hope we have helped the reader in grasping their potential for cryptographic aplications. While it is hard to predict how the field will evolve, we can for sure identify a number of interesting problems on the frontier between cryptography and group theory:  


\begin{itemize}
\item Babai's conjecture that short paths exist in Cayley graphs of finite simple groups is a widely-studied open problem in group theory. The Factorization Problem, equivalent to finding preimages for Cayley hash functions, requires \emph{constructing} such short paths in Cayley graphs. For cryptographic applications it is desirable to either find a situation in which the Factorization Problem is computationally infeasible, or to show that it is always feasible, as discussed in Section~\ref{sect: Cayley proofs}. Progress in constructing short enough paths would imply progress on Babai's conjecture.

\item Logarithmic signatures are a possible source of useful trapdoor functions for public-key cryptography, but there is more work to be done on understanding and constructing them. One direction, discussed in Section~\ref{sect: LS constructions}, is to find an algorithm that can produce wild logarithmic signatures, especially one which can also provide a rewriting in terms of tame ones. Another, discussed in Section~\ref{sect: MLS proofs}, is to determine whether all finite groups have minimal length logarithmic signatures. This question has been reduced to simple groups, and the MLS Conjecture that minimal length logarithmic signatures exist for all simple groups remains open in some cases.

\item Ostrovsky and Skeith \cite{OstSk2008} show how to convert a homomorphic encryption procedure on any finite simple group to a fully homomorphic encryption procedure on a ring by constructing $NAND$ in the finite simple groups. As discussed in Section~\ref{sect: FHE simple groups}, this opens up the question of finding secure homomorphic encryption on a finite simple group.

\item The Hidden Subgroup Problem is central to post-quantum cryptography. As discussed in Section~\ref{sect: HSP}, understanding the hardness of HSP for symmetric groups could be useful in the analysis of post-quantum group-based cryptographic primitives. 
\end{itemize}
\section*{Acknowledgement} We thank Greg Kuperberg for generously sharing his ideas on HSP on the non-abelian finite simple groups.

DK is grateful to MGV for the kind invitation to Madrid in February 2020 and great hospitality and generous support; when DK proposed the topic of this manuscript in that research trip.
\bibliographystyle{alpha}
\bibliography{BibSGCrypto}
\end{document}